\documentclass[12pt]{article}
\usepackage{amscd,amsfonts,amssymb,amsmath,latexsym,array,hhline}
\usepackage[dvips]{graphicx}

\usepackage{theorem}

\renewcommand{\bar}{\overline}
\renewcommand{\tilde}{\widetilde}
\renewcommand{\hat}{\widehat}

\renewcommand{\rho}{r}

\newcommand{\LL}{{\mathbb L}}

\newcommand{\NN}{{\mathbb N}}
\newcommand{\RR}{{\mathbb R}}

\renewcommand{\phi}{\varphi}

\numberwithin{equation}{section}

\renewcommand{\section}[1]{\refstepcounter{section}
                            {\noindent\large\bf\thesection. #1}}

\makeatletter

\renewcommand{\section}{\@startsection{section}{1}{0pt}{30pt}{6pt}{\large\bf}}
\def\dot{\hspace{-16pt}.\hspace{-2pt} }

\renewcommand{\@makefnmark}{}
\renewcommand{\@cite}[2]{[{#1\if@tempswa ; #2\fi}]}

\newenvironment{proof}{{\vskip\baselineskip\noindent\textbf{Proof}}}%
{\hspace*{.1pt}\hspace*{\fill}$\square$\vskip\baselineskip}
\newtheorem{theorem}{Theorem}[section]

\makeatother

\begin{document}

\title{Optimal stopping problems in diffusion-type \\ models with running maxima and drawdowns}

\author{Pavel V. Gapeev\footnote{London School of Economics,
Department of Mathematics, Houghton Street, London WC2A 2AE, United
Kingdom; e-mail: p.v.gapeev[n.rodosthenous]{\char'100}lse.ac.uk}
\and Neofytos Rodosthenous$^{*}$\footnote{Supported by the
scholarship of the Alexander Onassis Public Benefit Foundation.}}
\date{}
\maketitle


\begin{abstract}
We study optimal stopping problems related to the pricing of perpetual American options in an
extension of the Black-Merton-Scholes model in which the dividend and volatility rates of the
underlying risky asset depend on the running values of its maximum and maximum drawdown.
The optimal stopping times of exercise are shown to be the first times at which the price of the underlying asset
exits some regions restricted by certain boundaries depending on the running values of the associated maximum
and maximum drawdown processes. We obtain closed-form solutions to the equivalent free-boundary problems
for the value functions with smooth fit at the optimal stopping boundaries and normal reflection
at the edges of the state space of the resulting three-dimensional Markov process. We derive first-order
nonlinear ordinary differential equations for the optimal exercise boundaries of the perpetual American standard options.
\end{abstract}





\section{\dot Introduction}

       The main aim of this paper is to present closed-form solutions
       to the discounted optimal stopping problem of (\ref{4V5b}) for the initial process $X$
       in a diffusion-type model containing its running maximum process $S$ and its running
       maximum drawdown process $Y$ defined in (\ref{4X4})-(\ref{4S4}).
       This problem is related to the option pricing theory in mathematical finance,
       where the process $X$ can describe the price of a risky asset (e.g. a stock)
       on a financial market. The value of (\ref{4V5b}) can therefore be interpreted as the
       rational (or no-arbitrage) price of a standard perpetual American call or put option in a
       diffusion-type extension of the Black-Merton-Scholes model
       (see, e.g. Shiryaev \cite[Chapter~VIII; Section~2a]{FM}, Peskir and Shiryaev
       \cite[Chapter~VII; Section~25]{PSbook}, and Detemple \cite{Det},
       for an extensive overview of other related results in the area).

       Optimal stopping problems 
       for running maxima of some diffusion processes given linear costs were studied by Jacka \cite{Jmax},
       Dubins, Shepp, and Shiryaev \cite{DSS}, and Graversen and Peskir \cite{GP1}-\cite{GP2}
       among others, with the aim of determining the best constants in the corresponding maximal inequalities.
       A complete solution of a general version of the same problem was obtained in Peskir \cite{Pmax},
       by means of the established maximality principle, which is equivalent to the superharmonic
       characterization of the value function. Discounted optimal stopping problems for certain
       payoff functions depending on the running maxima of geometric Brownian motions were
       initiated by Shepp and Shiryaev \cite{SS1}-\cite{SS2} and then considered by
       Pedersen \cite{Jesper} and Guo and Shepp \cite{GuoShepp} among others,
       with the aim of computing rational values of perpetual American lookback (Russian) options.
       More recently, Guo and Zervos \cite{GuoZer} derived solutions for discounted optimal
       stopping problems related to the pricing of perpetual American options with certain payoff functions
       depending on the running values of both the initial diffusion process and its associated maximum.
       Glover, Hulley, and Peskir \cite{GHP} provided solutions of optimal stopping problems for integrals
       of functions depending on the running values of both the initial diffusion process and its associated minimum.
       The main feature of the resulting optimal stopping problems
       is that the normal-reflection condition holds for the value function at the diagonal of the state space
       of the two-dimensional continuous Markov process having the initial process and its running extremum
       as the components, which implies the characterization of the optimal boundaries as extremal solutions
       of one-dimensional first-order nonlinear ordinary differential equations.

       Asmussen, Avram, and Pistorius \cite{AAP} considered perpetual American options with payoffs depending
       on the running maximum of some L\'evy processes with two-sided jumps having phase-type distributions
       in both directions. Avram, Kyprianou, and Pistorius \cite{AKP} studied exit problems for spectrally
       negative L\'evy processes and applied the results to solving optimal stopping problems for payoff
       functions depending on the running values of the initial processes or their associated maxima.
       Optimal stopping games with payoff functions of such type were considered by Baurdoux and Kyprianou
       \cite{BK3}, within the framework of models based on spectrally negative L\'evy processes.
       Other complicated optimal stopping problems for the running maxima were considered in \cite{Gap}
       for a jump-diffusion model with compound Poisson processes with exponentially distributed jumps
       and by Ott \cite{Ott} (see also Ott \cite{Ott_thesis}) for a model based on spectrally negative L\'evy processes.
       More recently, Peskir \cite{Pe5a}-\cite{Pe5b} studied optimal stopping problems for three-dimensional Markov
       processes having the initial diffusion process as well as its maximum and minimum as the state space components.
       It was shown that the optimal boundary surfaces depending on the maximum and minimum of the initial process provide
       the maximal and minimal solutions of the associated systems of first-order non-linear partial differential equations.

In this paper, we obtain closed-form solutions to the perpetual American standard options pricing problem in an extension of the Black-Merton-Scholes model with path-dependent coefficients. The underlying asset price dynamics are described by a geometric diffusion-type process $X$ with local drift and diffusion coefficients, which essentially depend on the running values of the maximum process $S$ and the maximum drawdown process $Y$. It is shown that the optimal exercise times are the first times at which the process $X$ exits some regions restricted by certain boundaries depending on the running values of $S$ and $Y$. The process $Y$ represents the maximum of the difference between the running values of the underlying asset price and its maximum and can therefore be interpreted as the maximum of the market depth.
We derive closed-form expressions for the value functions as solutions of the equivalent free-boundary problems and apply the maximality principle from \cite{Pmax} to describe the optimal boundary surfaces as the {\it maximal} solutions of first-order nonlinear ordinary differential equations. The starting points for these surfaces at the edges of the three-dimensional state space of $(X, S, Y)$ are specified from the solutions of the corresponding optimal stopping problem for the two-dimensional Markov process $(X, S)$ in a corresponding model in which the coefficients of the process $X$ depend only on the running maximum process $S$.
The Laplace transforms of the drawdown process and other related characteristics associated with certain classes of the initial processes such as some diffusion models and spectrally positive and negative L\'evy processes, were studied by Pospisil, Vecer, and Hadjiliadis \cite{PVH} and by Mijatovic and Pistorius \cite{MP}, respectively.
Diffusion-type processes with given joint laws for the terminal level and supremum at an independent exponential time were
constructed in Forde \cite{Forde}, by allowing the diffusion coefficient to depend on the running value of the initial process
and its running minimum. Other important characteristics for such diffusion-type processes were recently derived by Forde, Pogudin, and Zhang \cite{FPZ}.

       The paper is organized as follows.
       In Section 2, we formulate the associated optimal stopping problems for a necessarily three-dimensional continuous Markov process, which has the underlying asset price and the running values of its maximum and maximum drawdown as the state space components.
       The resulting optimal stopping problem are reduced to their equivalent free-boundary problem for the value function, which satisfies the smooth-fit conditions at the stopping boundaries and the normal-reflection conditions at the edges of the state space of the three-dimensional process.
       In Section 3, we obtain closed-form solutions of the associated free-boundary problem in which the sought boundaries are found as appropriate systems of arithmetic equations or first-order nonlinear ordinary differential equations, where we specify the starting values for the latter at the edges of the three-dimensional state space.
       In Section 4, by applying the change-of-variable formula with local time on surfaces, we verify that the resulting solutions of the free-boundary problems provide the expressions for the value functions and the optimal stopping boundaries for the underlying asset price process in the initial problem. The main results of the paper are stated in Theorems \ref{4call}
       and \ref{4put}.


\section{\dot Preliminaries}

In this section, we introduce the setting and notation of the three-dimensional optimal stopping problems
which are related to the pricing of perpetual American standard options and formulate the equivalent free-boundary problems.

     \vspace{6pt}

     {\bf 2.1. Formulation of the problem.} For a precise formulation of the problem,
     let us consider a probability space $(\Omega, {\cal F}, P)$
     with a standard Brownian motion $B=(B_t)_{t \ge 0}$.
     Assume that there exists a process $X=(X_t)_{t \ge 0}$ given by
       \begin{equation}
       \label{4X4}
       X_t = x \, \exp \bigg( \int_0^t \Big( r - \delta(S_u, Y_u) -
       \frac{\sigma^2(S_u, Y_u)}{2} \Big) \, du + \int_0^t \sigma(S_u, Y_u) \, dB_u \bigg)
       \end{equation}
    where $\delta(s, y) > 0$ and $\sigma(s, y) > 0$
    are continuously differentiable bounded functions on $[0, \infty]^2$.
    It follows that the process $X$ solves the stochastic differential equation
       \begin{equation}
       \label{4dX4}
       dX_t = (r-\delta(S_t, Y_t)) \, X_t \, dt + \sigma(S_t, Y_t) \, X_t \, dB_t
       \quad (X_0=x)
       \end{equation}
     where $x > 0$ is given and fixed. Here, the associated with $X$
     {\it running maximum} process $S=(S_t)_{t \ge 0}$ and the corresponding
     {\it running maximum drawdown} process $Y=(Y_t)_{t \ge 0}$ are defined by
        \begin{equation}
        \label{4S4}
        S_t = s \vee \max_{0 \le u \le t} X_u \quad \text{and} \quad Y_t = y \vee \max_{0 \le u \le t} (S_u - X_u)
        \end{equation}
     for arbitrary $0 < s - y \le x \le s$, so that $X$ is a diffusion-type process representing
     a unique solution of the stochastic differential equation in (\ref{4dX4}) (see, e.g. \cite[Chapter~IV, Theorem~4.6]{LS}). 

     The main purpose of the present paper is
     to derive a closed-form solution to the optimal stopping
     problem for the time-homogeneous (strong) Markov
     process $(X, S, Y)=(X_t, S_t, Y_t)_{t \ge 0}$ given by
        \begin{equation}
        \label{4V5b}
        V_*(x, s, y) = \sup_{\tau} E_{x, s, y}
        \big[ e^{- r \tau} \, G(X_{\tau})
\big]
        \end{equation}
       for any $(x, s, y) \in E^3$, where the supremum is taken over all stopping times $\tau$
       with respect to the natural filtration of $X$,
       and the payoff function is either $G(x)=(L-x)^+$ or $G(x)=(x-K)^+$,
       for some given constants $0 < L, K < \infty$.
       Here, $E_{x, s, y}$ denotes the expectation
       under the assumption that
       the (three-dimensional) process $(X, S, Y)$ defined in
       (\ref{4X4})-(\ref{4S4}) starts at $(x, s, y) \in E^3$,
       where $E^3 = \{ (x, s, y) \in \RR^3 \, | \, 0 < s-y \leq x \leq s \}$ is
       the state space of the process $(X, S, Y)$.
       We assume that the process $X$ describes the price of a risky asset
       on a financial market, where $r$ is the riskless interest rate,
       $\delta(s, y)$ is the dividend rate paid to the asset holders,
       and $\sigma(s, y)$ is the volatility rate.
       The value of (\ref{4V5b}) is then actually a {\it rational} (or {\it no-arbitrage})
       price of a perpetual American put or call option with payoff function
       $G(x) = (L - x)^+$ or $G(x) = (x - K)^+$, where the expectation
       is taken under the (unique) martingale measure (see, e.g. \cite[Chapter~VII, Section~3g]{FM}).

       \vspace{6pt}

{\bf 2.2. The structure of the optimal stopping times.}
It follows from the general theory of optimal stopping problems
for Markov processes (see, e.g. \cite[Chapter~I, Section~2.2]{PSbook}) that the
optimal stopping time in the problem of (\ref{4V5b}) is given by
\begin{equation}
\label{4tau0} \tau_* = \inf \{t \ge 0 \, | \, V_*(X_t, S_t, Y_t) = G(X_t) \}
\end{equation}
so that, taking into account the structure of the payoff function
$G(x) = (L - x)^+$ or $G(x) = (x - K)^+$ in (\ref{4V5b}),
we further assume that the optimal stopping time from (\ref{4tau0}) takes either the form
\begin{equation}
\label{4tau}
\tau_* = \inf \{t \ge 0 \, | \, X_t \leq a_*(S_t, Y_t) \}
\quad \text{or} \quad
\tau_* = \inf \{t \ge 0 \, | \, X_t \geq b_*(S_t, Y_t) \}
\end{equation}
for some function $0 < a_*(s, y) < L$ or $K < b_*(s, y) < \infty$ to be determined.
This assumption means that the set
\begin{equation}
\label{4C'}
C' = \{(x,s,y) \in E^3 \, | \, a_*(s, y) < s-y  
\;\; \text{or} \;\; s < b_*(s, y) \}
\end{equation}
belongs to the continuation region for the optimal stopping problem of (\ref{4V5b})
which is given by
\begin{equation}
\label{4C}
C_* = \{(x,s,y) \in E^3 \, | \, a_*(s,y) < x 
\;\; \text{or} \;\; x < b_*(s,y) \}
\end{equation}
and the corresponding stopping region is the closure of the set
\begin{equation}
\label{4D}
D_* = \{(x,s,y) \in E^3 \, | \, x < a_*(s,y) 
\;\; \text{or} \;\; b_*(s,y) < x \}.
\end{equation}


\vspace{3pt}

{\bf 2.3. The free-boundary problem.} By means of standard arguments based on the application of
It\^o's formula, it is shown that the infinitesimal operator $\LL$ of the
process $(X, S, Y)$ acts on a function $F(x, s, y)$ from the class
$C^{2,1,1}$ on the interior of $E^3$ according to the rule
\begin{equation}
\label{4LF}
(\LL F)(x, s, y) = (r - \delta(s, y)) \, x \, \partial_x F(x, s, y)
+ \frac{\sigma^2(s, y)}{2} \, x^2 \, \partial^2_{xx} F(x, s, y)
\end{equation}
for all $0 < s-y < x < s$.
It follows from
the fact that both payoff functions $G(x)=(L - x)^+$ and $G(x)=(x - K)^+$ are convex
that the value function $V_*(x, s, y)$ is convex in the variable $x$,
and thus, it is continuous in $x$ on the interval $(0, \infty)$.
      In order to find analytic expressions for the unknown
      value function $V_*(x, s, y)$ from (\ref{4V5b}) and the unknown
      boundary $a_*(s, y)$ or $b_*(s, y)$ from (\ref{4tau}),
      let us build on the results of general theory of optimal stopping
      problems for Markov processes (see, e.g. 
      \cite[Chapter~IV, Section~8]{PSbook}). We can reduce the optimal
      stopping problem of (\ref{4V5b}) to the equivalent free-boundary problem
      for $V_*(x, s, y)$ with $a_*(s, y)$ or $b_*(s, y)$ given by
       \begin{align}
       \label{4LV}
       &({\LL} V)(x, s, y) = r \, V(x, s, y)
       \quad \text{for} \quad (x, s, y) \in C \\
       \label{4CF}
       &V(x, s, y) \big|_{x = a(s, y)+} = L - a(s, y) \quad \text{or} \quad
       V(x, s, y) \big|_{x = b(s, y)-} = b(s, y) - K \\
       \label{4VD}
       &V(x, s, y) = (L - x)^+ \quad \text{or} \quad V(x, s, y) = (x - K)^+ \quad \text{for} \quad (x, s, y) \in D \\
       \label{4VC}
       &V(x, s, y) > (L - x)^+ \quad \text{or} \quad V(x, s, y) > (x - K)^+ \quad \text{for} \quad (x, s, y) \in C \\
       \label{4LVD}
       &({\LL} V)(x, s, y) < r \, V(x, s, y) \quad \text{for} \quad (x, s, y) \in D
      \intertext{where $C$ and $D$ are defined as $C_*$ and $D_*$ in (\ref{4C}) and (\ref{4D})
      with $a(s, y)$ and $b(s, y)$ instead of $a_*(s, y)$ and $b_*(s, y)$, respectively,
      and the instantaneous-stopping conditions in (\ref{4CF}) are satisfied,
      when $s - y \le a(s, y)$ or $b(s, y) \le s$, respectively, for each $0 < y < s$.
      Observe that the superharmonic characterization of the value function (see \cite{Dyn} 
	  and \cite[Chapter~IV, Section~9]{PSbook}) implies that $V_*(x, s, y)$ is the smallest function
      satisfying (\ref{4LV})-(\ref{4VC}), with the boundary $a_*(s,y)$ or $b_*(s, y)$.
      Moreover, we further assume
for the left-hand system of (\ref{4LV})-(\ref{4LVD}),
that the smooth-fit and normal-reflection conditions}
       \label{4NRS}
       &\partial_x V(x, s, y) \big|_{x = a(s, y)+} = - 1 \quad
       \text{and} \quad \partial_s V(x, s, y) \big|_{x = s-} = 0 \\
       \intertext{hold, when $s - y \leq a(s, y) < s$, 
       and for the right-hand system of (\ref{4LV})-(\ref{4LVD}), that the normal-reflection and smooth-fit conditions}
       \label{4NRY}
       &\partial_y V(x, s, y) \big|_{x = (s-y)+} = 0 \quad \text{and} \quad
       \partial_x V(x, s, y) \big|_{x = b(s, y)-} = 1
       \end{align}
       hold, when 
$s - y < b(s, y) \leq s$, for each $0 < y < s$.
       Otherwise,
the normal-reflection
       conditions in the right-hand part of (\ref{4NRS}) and in the left-hand part
       of (\ref{4NRY}) hold, when either $a(s, y) < s - y$ or $s < b(s, y)$, respectively, for each $0 < y < s$.
%

Note that, when $\delta(s, y) = \delta(s)$ and
$\sigma(s, y) = \sigma(s)$ holds in (\ref{4X4})-(\ref{4dX4}),
the value function $V_*(x, s, y) = U_*(x, s)$ 
together with the boundary $a_*(s, y) = g_*(s)$ or $b_*(s, y) = h_*(s)$,
satisfy the left-hand or the right-hand system of (\ref{4LV})-(\ref{4LVD}),
respectively.
Moreover, the smooth-fit and normal-reflection conditions
\begin{align}
\label{4NRA}
&\partial_x U(x, s) \big|_{x = g(s)+} = - 1 \quad
\text{and} \quad \partial_s U(x, s) \big|_{x = s-} = 0 \\
\intertext{hold in addition to the left-hand system of (\ref{4LV})-(\ref{4LVD}), when $0 < g(s) < s$, 
and the natural-boundary
and smooth-fit conditions}
\label{4NRB}
&U(x, s) \big|_{x = 0+} = 0 
\quad \text{and} \quad
\partial_x U(x, s) \big|_{x = h(s)-} = 1
\end{align}
hold in addition to the right-hand system of (\ref{4LV})-(\ref{4LVD}), when 
$K < h(s) \le s$, for each $s > 0$.
Otherwise,
the normal-reflection and natural-boundary conditions in the right-hand part of (\ref{4NRA}) and in the left-hand part of (\ref{4NRB}) hold, respectively for the latter system, when 
$s < h(s)$, for each $s > 0$.

\section{\dot Solution of the free-boundary problem}
  \setcounter{equation}{0}

\label{ch4.2}

In this section, we obtain closed-form expressions for the value functions $V_*(x, s, y)$
in (\ref{4V5b}) for the payoffs of standard put and call options, and derive
explicit expressions and first-order nonlinear ordinary differential equations for the optimal exercise boundaries
$a_*(s, y)$ and $b_*(s, y)$ from (\ref{4tau}), as solutions to the free-boundary problems
(\ref{4LV})-(\ref{4LVD}) with (\ref{4NRS}) and (\ref{4NRY}), respectively.

\vspace{6pt}

{\bf 3.1. The general solution of the free-boundary problem.}
We first observe that the general solution of the equation in (\ref{4LV}) has the form
\begin{equation}
\label{4V0}
V(x, s, y) = C_1(s, y) \, x^{\gamma_1(s, y)} + C_2(s, y) \, x^{\gamma_2(s, y)}
\end{equation}
where $C_i(s, y)$, $i = 1, 2$, are some arbitrary
continuously differentiable functions and $\gamma_2(s, y) < 0 < 1 < \gamma_1(s, y)$ are given by
\begin{equation}
\label{4gamma12}
\gamma_{i}(s, y) = \frac{1}{2} - \frac{r - \delta(s, y)}{\sigma^2(s, y)}
- (-1)^i \sqrt{ \bigg( \frac{1}{2} - \frac{r - \delta(s, y)}
{\sigma^2(s, y)} \bigg)^2 + \frac{2 r}{\sigma^2(s, y)}}
\end{equation}
for all $0 < y < s$.
Hence, applying the instantaneous-stopping conditions from (\ref{4CF}) to the function
in (\ref{4V0}), we get that the equalities
\begin{align}
\label{4B31a}
&C_1(s, y) \, a^{\gamma_1(s, y)}(s, y)
+ C_2(s, y) \, a^{\gamma_2(s, y)}(s, y) = L - a(s, y) \\
\label{4B31b}
&C_1(s, y) \, b^{\gamma_1(s, y)}(s, y)
+ C_2(s, y) \, b^{\gamma_2(s, y)}(s, y) = b(s, y) - K
\end{align}
hold, when $s-y \le a(s, y)$ and $b(s, y) \le s$, respectively,
for each $0 < y < s$. Moreover, using the smooth-fit conditions from the
left-hand part of (\ref{4NRS}) and the right-hand part of (\ref{4NRY}), we obtain
that the equalities
\begin{align}
\label{4B31aa}
&C_1(s, y) \, \gamma_1(s, y) \, a^{\gamma_1(s, y)}(s, y)
+ C_2(s, y) \, \gamma_2(s, y) \, a^{\gamma_2(s, y)}(s, y) = - a(s, y) \\
\label{4B31bb}
&C_1(s, y) \, \gamma_1(s, y) \, b^{\gamma_1(s, y)}(s, y)
+ C_2(s, y) \, \gamma_2(s, y) \, b^{\gamma_2(s, y)}(s, y) = b(s, y)
\end{align}
hold, when $s-y \le a(s, y)$ and $b(s, y) \le s$, respectively, for each $0 < y < s$.
Finally, applying the normal-reflection conditions from the right-hand part of (\ref{4NRS})
and the left-hand part of (\ref{4NRY}) to the function in (\ref{4V0}), we have that the equalities
\begin{align}
\label{4B31c}
&\sum_{i=1}^{2} \Big( \partial_s C_i(s, y) \, s^{\gamma_i(s,y)} + C_i(s,y)
\, \partial_s \gamma_i(s,y) \, s^{\gamma_i(s,y)} \ln s \Big) = 0 \\
\label{4B31d}
&\sum_{i=1}^{2} \Big( \partial_y C_i(s, y) \, (s-y)^{\gamma_i(s,y)} +
C_i(s, y) \, \partial_y \gamma_i(s,y) \, (s-y)^{\gamma_i(s,y)} \ln (s-y) \Big) = 0
\end{align}
hold, when either $a(s, y) < s$ or $s < b(s, y)$ and either $a(s, y) < s-y$ or $s-y < b(s, y)$,
respectively, for each $0 < y < s$.
Here, the partial derivatives $\partial_s \gamma_{i}(s, y)$ and $\partial_y \gamma_{i}(s, y)$
take the form
\begin{align}
\label{4gammas}
&\partial_s \gamma_{i}(s, y) = \phi(s, y) - (-1)^i \, \frac{
\phi(s, y) (\gamma_1(s, y) - \gamma_2(s, y)) \sigma^3(s, y) - 2 r \partial_s \sigma(s, y)}
{\sigma^2(s, y) \sqrt{(\gamma_1(s, y) - \gamma_2(s, y))^2 \sigma^2(s, y) + 2 r}} \\
\label{4gammay}
&\partial_y \gamma_{i}(s, y) = \psi(s, y) - (-1)^i \, \frac{
\psi(s,y) (\gamma_1(s, y) - \gamma_2(s, y)) \sigma^3(s, y) - 2 r \partial_y \sigma(s, y)}
{\sigma^2(s, y) \sqrt{(\gamma_1(s, y) - \gamma_2(s, y))^2 \sigma^2(s, y) + 2 r}}
\end{align}
for $i = 1, 2$, and the functions $\phi(s, y)$ and $\psi(s, y)$ are defined by
\begin{align}
\label{4phi}
&\phi(s, y) = \frac{\sigma(s, y) \partial_s \delta(s, y) + 2
(r - \delta(s, y)) \partial_s \sigma(s, y)}{\sigma^3(s, y)} \\
\label{4psi}
&\psi(s, y) = \frac{\sigma(s, y) \partial_y \delta(s, y) + 2
(r - \delta(s, y)) \partial_y \sigma(s, y)}{\sigma^3(s, y)}
\end{align}
for $0 < y < s$.

\vspace{6pt}



{\bf 3.2. The solution to the problem in the two-dimensional $(X, S)$-setting.}
We begin with the case in which $\delta(s, y) = \delta(s)$ and
$\sigma(s, y) = \sigma(s)$ holds in (\ref{4X4})-(\ref{4dX4}), and thus, we can define
the functions $\beta_i(s) = \gamma_i(s, y)$, $i = 1, 2$, as in (\ref{4gamma12}).
Then, the general solution $V(x, s, y) = U(x, s)$ of the equation in (\ref{4LV})
has the form of (\ref{4V0}) with $C_i(s, y) = D_i(s)$ and $\gamma_i(s, y) = \beta_i(s)$,
for $i = 1, 2$, and the stopping time takes the form of (\ref{4tau}) with either
the boundary $a_*(s, y) = g_*(s)$ or $b_*(s, y) = h_*(s)$, respectively. We further denote the state space
of the two-dimensional (strong) Markov process $(X, S)$ by $E^2 = \{ (x, s) \in \RR^2 \,
| \, 0 < x \leq s \}$ and its diagonal by $d^2 = \{ (x, s) \in \RR^2 \, | \, 0 < x = s \}$,
as well as recall that the second component of $(X, S)$ can only increase at $d^2$, that is,
when $X_t = S_t$ for $t \ge 0$.

\vspace{3pt}

{\bf (i) The case of call option.}
Let us first consider the payoff function $G(x) = (x-K)^+$ in (\ref{4V5b}).
In this case, taking into account the fact that $\beta_2(s) < 0 < 1 < \beta_1(s)$,
we observe that $D_2(s) = 0$ should hold in (\ref{4V0}), since otherwise $U(x, s) \to \pm \infty$
as $x \downarrow 0$, that must be excluded, by virtue of the obvious fact that the value function
in (\ref{4V5b}) is bounded at zero, due to the natural boundary condition in the left-hand part
of (\ref{4NRB}).
Hence, solving the system of equations in (\ref{4B31b}) and (\ref{4B31bb}) for the unknown function
$C_1(s, y) = D_1(s)$ with $C_2(s, y) = D_2(s) = 0$, we conclude that the function
$V(x, s, y) = U(x, s)$ in (\ref{4V0}) admits the representation
\begin{equation}
\label{4U33a}
U(x, s; h_*(s)) = \frac{h_*(s)}{\beta_1(s)} \, \Big( \frac{x}{h_*(s)} \Big)^{\beta_1(s)}
\quad \text{with} \quad h_*(s) = \frac{\beta_1(s) K}{\beta_1(s) - 1}
\end{equation}
for $0 < x < h_*(s) \leq s$, so that $h_*(s) > K \vee (r K/\delta(s))$ holds for all $s > K$.

In this case, we set ${\tilde s}_{0} = \infty$ and define a decreasing sequence
$({\tilde s}_n)_{n \in \NN}$ such that the boundary $h_*(s)$ from (\ref{4U33a})
exits the region $E^2$ at $({\tilde s}_{2l-1}, {\tilde s}_{2l-1})$ and enters
back to $E^2$ downwards at $({\tilde s}_{2l}, {\tilde s}_{2l})$. Namely, we define
${\tilde s}_{2l-1} = \sup \{ s < {\tilde s}_{2l-2} \, | \, h_*(s) > s \}$
and ${\tilde s}_{2l} = \sup \{ s < {\tilde s}_{2l-1} \, | \, h_*(s) \leq s \}$, $k \in \NN$,
whenever they exist, and put ${\tilde s}_{2l-1} = {\tilde s}_{2l} = 0$ otherwise.
Note that $K < {\tilde s}_{2l} < {\tilde s}_{2l-1} < \infty$, $l \in \NN$, by construction.
Then, the candidate value function admits the representation of (\ref{4U33a}) in the regions
\begin{equation} \label{4tilR2_2l-1}
{\tilde Q}^2_{2l-1} = \{ (x, s) \in E^2 \, | \, {\tilde s}_{2l-1} < s \le {\tilde s}_{2l-2} \}
\end{equation}
for $l \in \NN$.

On the other hand, the candidate value function $V(x, s, y) = U(x, s)$
takes the form of (\ref{4V0}) with $C_1(s, y) = D_1(s)$ solving the
first-order linear ordinary differential equation in (\ref{4B31c}) and
$C_2(s, y) = D_2(s) = 0$, in the regions
\begin{equation} \label{4tilR2_2l}
{\tilde R}^2_{2l} = \{(x, s) \in E^2 \, | \, {\tilde s}_{2l} < s \le {\tilde s}_{2l-1} \}
\end{equation}
for $l \in \NN$, which belong to $C'$ in (\ref{4C'}).
Note that the process $(X, S)$ can pass from the region ${\tilde R}^2_{2l}$
in (\ref{4tilR2_2l}) to the region ${\tilde R}^2_{2l-1}$ in (\ref{4tilR2_2l-1})
only through the point $({\tilde s}_{2l-1}, {\tilde s}_{2l-1})$, for $l \in \NN$. Thus, the candidate value function should be continuous at the point $({\tilde s}_{2l-1}, {\tilde s}_{2l-1})$, that is expressed by the equality
\begin{equation}
\label{4U33d}
D_1({\tilde s}_{2l-1}) \, ({\tilde s}_{2l-1})^{\beta_1({\tilde s}_{2l-1})}
= U({\tilde s}_{2l-1}+, {\tilde s}_{2l-1}+; h_*({\tilde s}_{2l-1}+))
\end{equation}
where the right-hand side is given by (\ref{4U33a}).
Hence, solving the first-order linear ordinary differential
equation in (\ref{4B31c}) for the unknown function $C_1(s, y) = D_1(s)$ with $C_2(s, y) = D_2(s) = 0$
and using the condition of (\ref{4U33d}), we obtain that the candidate value function $V(x, s, y) = U(x, s)$
in (\ref{4V0}) admits the expression
\begin{equation}
\label{4Uint2}
U(x, s; {\tilde s}_{2l-1}) = \exp \bigg( - \int_s^{{\tilde s}_{2l-1}} \beta_1'(q) \, \ln q \, dq \bigg)
\, \frac{({\tilde s}_{2l-1})^{1-\beta_1({\tilde s}_{2l-1})}}{\beta_1({\tilde s}_{2l-1})} \,
{x}^{\beta_1(s)}
\end{equation}
in the regions ${\tilde R}^2_{2l}$ given by (\ref{4tilR2_2l}), for $l \in \NN$.

\vspace{3pt}


{\bf (ii) The case of put option.}
Let us now consider the payoff function $G(x) = (L-x)^+$ in (\ref{4V5b}).
In this case, solving the system of equations in (\ref{4B31a}) and (\ref{4B31aa}) for the unknown functions
$C_i(s, y) = D_i(s)$,
$i = 1, 2$, 
we conclude that the function $V(x, s, y) = U(x, s)$ in (\ref{4V0}) admits the representation
\begin{equation}
\label{4U32a}
U(x, s; g_*(s)) = D_1(s; g_*(s)) \, x^{\beta_1(s)} + D_2(s; g_*(s)) \, x^{\beta_2(s)}
\end{equation}
for $0 < g_*(s) < x \leq s$, with
\begin{equation}
\label{4Di32a}
D_i(s; g_*(s)) = \frac{(\beta_{3-i}(s) - 1) g_*(s) - \beta_{3-i}(s) L}
{(\beta_i(s) - \beta_{3-i}(s)) g_*(s)^{\beta_i(s)}}
\end{equation}
for all $s > 0$ and $i = 1, 2$.
Hence, assuming that the boundary function $g_*(s)$ is continuously differentiable, we apply the
condition of (\ref{4B31c}) to the functions $C_i(s, y) = D_i(s; g_*(s))$, $i = 1, 2$, in (\ref{4Di32a})
and obtain that $g_*(s)$ solves the first-order nonlinear ordinary differential equation
\begin{align}
\label{4g'32a}
g'(s) &= \sum_{i=1}^{2}
\frac{((\beta_{3-i}(s)-1) g(s) - \beta_{3-i}(s)) g(s)}
{(\beta_i(s)-1) (\beta_{3-i}(s)-1) g(s) - \beta_i(s) \beta_{3-i}(s) L} \\
\notag
&\phantom{= \sum_{i=1}^{2} \;\:} \times \bigg( \frac{1}{\beta_{3-i}(s)
- \beta_i(s)} + \frac{(s/g(s))^{\beta_i(s)} \ln(s/g(s))}
{(s/g(s))^{\beta_{i}(s)}-(s/g(s))^{\beta_{3-i}(s)}} \bigg) \, \beta_i'(s)
\end{align}
where the derivatives $\beta_i'(s) = \partial_s \gamma_i(s, y)$,
$i = 1, 2$, are given by (\ref{4gammas}) with (\ref{4phi}).
Taking into account the fact that $\beta_i(s)$, $i = 1, 2$, and the boundary $g_*(s)$
are continuously differentiable functions in the neighborhood of infinity,
we observe that the function in (\ref{4U32a}) should satisfy the property
$U(x, s; g_*(s)) \to U(x, \infty; g_*(\infty))$ as $s \to \infty$,
for each $x > g_*(s)$. Thus, using the fact that $\beta_2(s) < 0 < 1 < \beta_1(s)$,
we obtain the expressions
\begin{equation}
\label{4U32b}
U(x, \infty; g_*(\infty)) = \frac{g_*(\infty)}{\beta_2(\infty)} \,
\Big( \frac{x}{g_*(\infty)} \Big)^{\beta_2(\infty)} \quad \text{and} \quad
g_*(\infty) = \frac{\beta_2(\infty) L}{\beta_2(\infty) - 1}
\end{equation}
for $x > g_*(\infty)$. The form of the function $U(x, \infty; g_*(\infty))$
and the boundary $g_*(\infty)$ in (\ref{4U32b}) follows from the fact that
$U(x, \infty; g_*(\infty)) \to \pm \infty$ should not hold as $x \to \infty$,
since the value function in (\ref{4V5b}) is bounded at infinity.
Observe that the expressions in (\ref{4U32b}) coincide
with the ones of the value function in the corresponding continuation region and
the exercise boundary of the perpetual American put option in the Black-Merton-Scholes
model with constant coefficients (see, e.g. \cite[Chapter~VIII, Section~2a]{FM}).

Let us now consider the maximal solution $g_*(s)$ of the first-order ordinary differential
equation in (\ref{4g'32a}) with the starting value $g_*(\infty)$ from (\ref{4U32b}) as $s \uparrow \infty$
and such that this solution stays below the curve $x = L \wedge (r L/\delta(s))$. 
Then, we put ${\hat s}_0 = \infty$ and define a decreasing sequence $({\hat s}_n)_{n \in \NN}$
such that the solution $g_*(s)$ of the equation in (\ref{4g'32a})
exits the region $E^2$ at the points $({\hat s}_{2k-1}, {\hat s}_{2k-1})$ and
enters $E^2$ downwards at the points $({\hat s}_{2k}, {\hat s}_{2k})$.
Namely, we define ${\hat s}_{2k-1} = \sup \{ s \le {\hat s}_{2k-2} \, | \, g_*(s) > s \}$
and ${\hat s}_{2k} = \sup \{ s \le {\hat s}_{2k-1} \, | \, g_*(s) \le s \}$,
$k \in \NN$, whenever they exist, and put ${\hat s}_{2k} = {\hat s}_{2k-1} = 0$
otherwise. Note that $0 < {\hat s}_{2k} < {\hat s}_{2k-1} < L$, $k \in \NN$, by construction.
Then, the candidate value function takes the form of (\ref{4U32a})-(\ref{4Di32a}) in the regions
\begin{equation} \label{4hatR2_2k-1}
{\hat Q}^2_{2k-1} = \{ (x, s) \in E^2 \, | \, {\hat s}_{2k-1} < s \le {\hat s}_{2k-2} \}
\end{equation}
for $k \in \NN$, and the boundary function $g_*(s)$ provides the maximal
solution of the equation in (\ref{4g'32a}) started at $g_*(\infty)$ from (\ref{4U32b})
and such that it stays strictly below the curve $x = L \wedge (r L/\delta(s))$. 
Finally, we note that the candidate value function should be given by the condition
in the left-hand part of (\ref{4VD}) in the regions
\begin{equation} \label{4hatR2_2k}
{\hat Q}^2_{2k} = \{ (x, s) \in E^2 \, | \, {\hat s}_{2k} < s \le {\hat s}_{2k-1} \}
\end{equation}
for $k \in \NN$, which belong to the stopping region $D_*$ in (\ref{4D}).

\vspace{6pt}


{\bf 3.3. The solution to the problem in the three-dimensional $(X,S,Y)$-setting.}
We now continue with the general form of the coefficients $\delta(s, y)$ and $\sigma(s, y)$
in (\ref{4X4})-(\ref{4dX4}), and thus, of the functions $\gamma_i(s, y)$, $i = 1, 2$, from (\ref{4gamma12}).
We denote the border planes of the state space $E^3$ by $d^3_1 = \{(x, s, y) \in \RR^3 \, | \, 0 < x = s \}$
and $d^3_2 = \{(x, s, y) \in \RR^3 \, | \, 0 < x = s-y \}$, as well as recall that the second and third components
of the process $(X, S, Y)$ can increase only at the planes $d^3_1$ and $d^3_2$, that is, when $X_t = S_t$ and
$X_t = S_t - Y_t$ for $t \ge 0$, respectively.

\vspace{3pt}


{\bf (i) The case of call option.}
Let us now get back to the payoff function $G(x) = (x - K)^+$ in (\ref{4V5b}). In this case, solving
the system of equations in (\ref{4B31b}) and (\ref{4B31bb}), we conclude
that the function in (\ref{4V0}) admits the representation
\begin{equation}
\label{4V33}
V(x, s, y; b_*(s, y)) = C_1(s, y; b_*(s, y)) \, x^{\gamma_1(s, y)} +
C_2(s, y; b_*(s, y)) \, x^{\gamma_2(s, y)}
\end{equation}
for $0 < s-y \le x < b_*(s, y) \leq s$ and $s>K$, with
\begin{equation}
\label{4Ci33}
C_i(s, y; b_*(s, y)) = \frac{(\gamma_{3-i}(s, y) - 1) b_*(s, y) - \gamma_{3-i}(s, y) K}
{(\gamma_{3-i}(s, y) - \gamma_i(s, y)) b_*(s, y)^{\gamma_i(s, y)}}
\end{equation}
for all $0 < y < s$ and $i = 1, 2$.
Hence, assuming that the boundary function $b_*(s, y)$ is continuously differentiable,
we apply the condition of (\ref{4B31d}) to the functions $C_i(s, y) = C_i(s, y; b_*(s, y))$, $i = 1, 2$,
in (\ref{4Ci33}) to obtain that $b_*(s, y)$ solves the first-order nonlinear ordinary differential equation
\begin{align}
\label{4g'33}
&\partial_y b(s, y) = \sum_{i=1}^{2}
\frac{((\gamma_{3-i}(s, y) - 1) b(s, y) - \gamma_{3-i}(s, y) K) b(s, y)}
{(\gamma_i(s, y) - 1) (\gamma_{3-i}(s, y) - 1) b(s, y) - \gamma_i(s, y) \gamma_{3-i}(s, y) K} \\
&\times \bigg( \frac{1}{\gamma_{3-i}(s, y) - \gamma_i(s, y)} +
\frac{((s - y)/b(s,y))^{\gamma_{i}(s,y)} \, \ln {((s - y)/b(s,y))}}
{((s - y)/b(s,y))^{\gamma_i(s, y)} - ((s - y)/b(s,y))^{\gamma_{3-i}(s, y)}} \bigg)
\, \partial_y \gamma_{i}(s, y) \nonumber
\end{align}
for $0 < y < s$, where the partial derivatives $\partial_y \gamma_{i}(s, y)$,
$i = 1, 2$, are given by (\ref{4gammay}) with (\ref{4psi}).

Since the functions $\delta(s, y)$ and $\sigma(s, y)$ are assumed
to be continuously differentiable and bounded, it follows that
the limits $\delta(s, s-)$ and $\sigma(s, s-)$ exist for each $s > 0$.
Then, the limits $\gamma_i(s, s-)$ can be identified with the functions
$\beta_i(s)$, $i = 1, 2$, from Subsection 3.2 above, and the function
in (\ref{4V33}) should satisfy the property $V(x, s, y; b_*(s, y)) \rightarrow V(x, s, s-; b_*(s, s-))$
as $y \uparrow s$, for each $s-y \leq x < b_*(s,y)$. Thus, taking into account
the fact that $\gamma_2(s, y) < 0 < 1 < \gamma_1(s, y)$, we conclude that
the equalities
\begin{equation}
\label{4Vss}
V(x, s, s-; b_*(s, s-)) = U(x, s; b_*(s, s-)) \quad \text{and} \quad
b_*(s, s-) = h_*(s)
\end{equation}
hold for $0 < x < b_*(s, s-)$ and $s > K$, with $U(x, s; h_*(s))$ and $h_*(s)$
given by (\ref{4U33a}), since otherwise
$V(x, s, s-; b_*(s, s-)) \to \pm \infty$ as $x \downarrow 0$, that must be
excluded by virtue of the obvious fact that the value function in (\ref{4V5b})
is bounded at zero.

For any $s > K$ fixed, let us now consider the solution $b_*(s, y)$ of
(\ref{4g'33}) started from the value $h_*(s)$ given by (\ref{4U33a}) at
$y \uparrow s$, given that this solution stays strictly above the
surface $x = K \vee (r K/\delta(s, y))$.
Then, we put ${\tilde y}_0(s) = s$ and define a decreasing sequence
$({\tilde y}_n(s))_{n \in \NN}$ such that
${\tilde y}_{2l-1}(s) = \sup \{ y < {\tilde y}_{2l-2}(s) \, | \, b_*(s, y) > s \}$
and ${\tilde y}_{2l}(s) = \sup \{ y < {\tilde y}_{2l-1}(s) \, | \, b_*(s, y) \le s \}$,
whenever they exist, and put ${\tilde y}_{2l-1}(s) = {\tilde y}_{2l}(s) = 0$, $l \in \NN$, otherwise.
Moreover, we can also define a decreasing sequence $({\hat y}_n(s))_{n \in \NN}$ such
that the boundary $b_*(s,y)$ exits the region $E^3$ from the side of $d^3_2$ at the points
$(s-{\hat y}_{2k-1}(s), s, {\hat y}_{2k-1}(s))$ and
enters $E^3$ downwards at the points $(s-{\hat y}_{2k}(s), s, {\hat y}_{2k}(s))$.
Namely, we put ${\hat y}_0(s) = s$ and define
${\hat y}_{2k-1}(s) = \sup \{ y < {\hat y}_{2k-2}(s) \, | \, b_*(s, y) < s - y \}$
and ${\hat y}_{2k}(s) = \sup \{ y < {\hat y}_{2k-1}(s) \, | \, b_*(s, y) \ge s - y \}$,
whenever such points exist, and put ${\hat y}_{2k-1}(s) = {\hat y}_{2k}(s) = 0$
otherwise, for $k \in \NN$. Note that $0 < {\hat y}_{2k}(s) < {\hat y}_{2k-1}(s) < s - K$,
$k \in \NN$, by construction.
Therefore, the candidate value function admits the expression in (\ref{4V33})-(\ref{4Ci33})
in either the region
\begin{equation}
\label{4tilR3_2l-1}
{\tilde R}^3_{2l-1} = \{ (x,s,y) \in E^3 \, | \, {\tilde y}_{2l-1}(s) < y \leq \min_{k \in \NN} \{ {\hat y}_{2k-2}(s) \, | \, {\tilde y}_{2l-1}(s) < {\hat y}_{2k-2}(s) \} \wedge {\tilde y}_{2l-2}(s) \}
\end{equation}
or
\begin{equation}
\label{4hatR3_2k-1}
{\hat R}^3_{2k-1} = \{ (x, s, y) \in E^3 \, | \, {\hat y}_{2k-1}(s) < y \leq  \min_{l \in \NN} \{ {\tilde y}_{2l-1}(s) \, | \, {\tilde y}_{2l-1}(s) < {\hat y}_{2k-2}(s) \} \wedge {\hat y}_{2k-2}(s) \}
\end{equation}
for $k,l \in \NN$, and the boundary $b_*(s, y)$ provides the unique solution of the equation in (\ref{4g'33})
started from the value $b_*(s, s-) = h_*(s)$ from (\ref{4U33a}) and such that this solution stays strictly above the
surface $x = K \vee (r K/\delta(s, y))$ 
(see Figure 1 below).

\begin{picture}(160,110)
\put(20,15){\begin{picture}(120,95)

\put(0,0){\line(1,0){120}} \put(120,0){\line(0,1){95}}
\put(0,0){\line(0,1){95}} \put(0,95){\line(1,0){120}}   

\put(10,10){\vector(1,0){100}}
\put(10,10){\vector(0,1){80}}    


\put(81,9){\line(0,1){72}} 

\put(80,6){$s$}

\put(9,81){\line(1,0){72}} 

\put(10,81){\line(1,-1){71}} 

\put(6,80){${s}$}

\put(70,40){\line(-1,0){19}}
\put(70,39.9){\line(-1,0){19}}
\put(70,40.1){\line(-1,0){19}}

\put(51,40){\line(-1,1){10}}
\put(51,39.9){\line(-1,1){10}}
\put(51,40.1){\line(-1,1){10}}
\put(51,39.8){\line(-1,1){10}}
\put(51,40.2){\line(-1,1){10}}

\put(75,43){\line(-1,0){27}}
\put(75,42.9){\line(-1,0){27}}
\put(75,43.1){\line(-1,0){27}}

\put(60,46.5){\line(-1,0){15.5}}
\put(60,46.4){\line(-1,0){15.5}}
\put(60,46.6){\line(-1,0){15.5}}

\put(81,50){\line(-1,0){40}}
\put(81,49.9){\line(-1,0){40}}
\put(81,50.1){\line(-1,0){40}}

\qbezier[5](81,50)(83,50)(85,50)
\qbezier[80](85,10)(85,47.5)(85,85)
\qbezier[75](10,85)(47.5,85)(85,85)
\qbezier[95](10,85)(47.5,50)(85,10)

\put(85,9){\line(0,1){2}} \put(84,6){$s'$}
\put(9,85){\line(1,0){2}} \put(6,84){$s'$}

\put(29,20){\vector(1,1){20.5}} \put(19,17){$x = s-y$}

\put(63,40){${\bf .}$}
\put(63,39.9){${\bf .}$}	\put(63.1,39.9){${\bf .}$}	\put(62.9,39.9){${\bf .}$}
\put(63,39.8){${\bf .}$}	\put(63.1,39.8){${\bf .}$} 	\put(62.9,39.8){${\bf .}$}
\put(63,39.7){${\bf .}$}	\put(63.1,39.7){${\bf .}$}	\put(62.9,39.7){${\bf .}$}
\put(63,39.6){${\bf .}$}

\put(61,36){$(x,s,y)$}

\put(0,0){\qbezier(69,81)(71,65)(81,60)} 
\put(0,0){\qbezier(81,60)(99,49)(81,35)}
\put(0,0){\qbezier(81,35)(71,28)(59,10)}

\qbezier[130](49,81)(79,49)(39,10) 

\qbezier[71](10,60)(45.5,60)(81,60)      
\qbezier[71](10,34.8)(45.5,34.8)(81,34.8)
\qbezier[58](10,22.7)(39,22.7)(68,22.7)

\put(.5,59){${\tilde y}_{1}(s)$}
\put(.5,33.8){${\tilde y}_{2}(s)$}
\put(.5,21.7){${\hat y}_{1}(s)$}

\qbezier[71](69,81)(69,45.5)(69,10)
\put(69,9){\line(0,1){2}} \put(64,6){$h_*(s)$}

\put(95,79){\vector(-2,-1){21}} \put(96,78){$b_*(s,y)$}

\put(45,71){\vector(2,1){8}} \put(36,68){$\frac{rK}{\delta(s,y)}$}
\put(35,9){\line(0,1){2}} \put(32.5,6){$K$}

\put(6,88){$y$} \put(108,6){$x$}
\end{picture}}
\put(17,10){\small{{\bf Figure 1.} A computer drawing of the
state space of the process $(X,S,Y)$,}} \\
\put(17,5) {\small{for some $s$ fixed, which increases to
$s'$, and the boundary function $b_*(s,y)$.}}
\end{picture}


On the other hand, the candidate value function
takes the form of (\ref{4V0}) with $C_{i}(s, y)$, $i = 1, 2$,
solving the linear system of first-order partial differential equations in (\ref{4B31c}) and (\ref{4B31d}),
in the regions
\begin{equation}
\label{4tilR3_2l}
{\tilde R}^3_{2l} = \{ (x,s,y) \in E^3 \, | \, {\tilde y}_{2l}(s) < y \leq {\tilde y}_{2l-1}(s) \}
\end{equation}
for $l \in \NN$, which belong to $C'$ in (\ref{4C'}).
Note that, the process $(X, S, Y)$ can enter the region ${\tilde R}^3_{2l}$
in (\ref{4tilR3_2l}) from one of the regions ${\tilde R}^3_{2l+1}$ in (\ref{4tilR3_2l-1})
or ${\hat R}^3_{2k-1}$ in (\ref{4hatR3_2k-1}), for some $k \in \NN$,
only through the point $(s - {\tilde y}_{2l}(s), s, {\tilde y}_{2l}(s))$ and can exit
the region ${\tilde R}^3_{2l}$ passing to the region ${\tilde R}^3_{2l-1}$ only through the
point $(s - {\tilde y}_{2l-1}(s), s, {\tilde y}_{2l-1}(s))$, by hitting the plane $d^3_2$,
so that increasing its third component $Y$.
Thus, the candidate function should be continuous at the points $(s - {\tilde y}_{2l}(s), s, {\tilde y}_{2l}(s))$
and $(s - {\tilde y}_{2l-1}(s), s, {\tilde y}_{2l-1}(s))$, that is expressed by the equalities
\begin{align}
\label{4condcall1}
&C_1(s, {\tilde y}_{2l}(s)+) \, ((s - {\tilde y}_{2l}(s))-)^{\gamma_1(s, {\tilde y}_{2l}(s)+)}
+ C_2(s, {\tilde y}_{2l}(s)+) \, ((s - {\tilde y}_{2l}(s))-)^{\gamma_2(s, {\tilde y}_{2l}(s)+)} \\ \nonumber
&= V(s - {\tilde y}_{2l}(s), s, {\tilde y}_{2l}(s); b(s, {\tilde y}_{2l}(s))) \\
\label{4condcall3}
&C_1(s, {\tilde y}_{2l-1}(s)) \, (s - {\tilde y}_{2l-1}(s))^{\gamma_1(s, {\tilde y}_{2l-1}(s))} +
C_2(s, {\tilde y}_{2l-1}(s)) \, (s - {\tilde y}_{2l-1}(s))^{\gamma_2(s, {\tilde y}_{2l-1}(s))} \\ \nonumber
&= V((s - {\tilde y}_{2l-1}(s))-, s, {\tilde y}_{2l-1}(s)+; b_*(s, {\tilde y}_{2l-1}(s)+)) 
\end{align}
for $s > K$ and $l \in \NN$, where the right-hand sides are given by (\ref{4V33})-(\ref{4Ci33})
with $b_*(s, {\tilde y}_{2l-1}(s)+) = b_*(s, {\tilde y}_{2l}(s)) = s$.
However, if $b_*(s, s-) = h_*(s) > s$ holds with $h_*(s)$ given by (\ref{4U33a}),
then we have ${\tilde y}_{1}(s) = s-$ and the condition of (\ref{4condcall3}) for $l = 1$, changes
its form to $C_2(s, s-) = 0$ for $s > K$,
since otherwise $V(x, s, s-) \rightarrow \pm \infty$ as $x \downarrow 0$, that must be excluded
by virtue of the obvious fact that the value function in (\ref{4V5b}) is bounded at zero.

In addition, the process $(X, S, Y)$ can exit the region ${\tilde R}^3_{2l}$ in (\ref{4tilR3_2l})
passing to the stopping region $D_*$ from (\ref{4D}) only through the point
$({\bar s}(y), {\bar s}(y), y)$, by hitting the plane $d^3_1$, so that increasing its second
component $S$ until it reaches the value ${\bar s}(y) = \inf \{ q > s \, | \, b_*(q, y) \leq q \}$.
Observe that the boundary $b_*(q, y)$ provides the unique solution of the equation in (\ref{4g'33})
with the starting value $b_*(q, q-) = h_*(q)$, for each $q \le {\bar s}(y)$, given that
this solution stays strictly above the surface $x = K \vee (r K/\delta(q, y))$.
Then, the candidate value function should be continuous at the point $({\bar s}(y), {\bar s}(y), y)$,
that is expressed by the equality
\begin{align}
\label{4condcall2}
&C_1({\bar s}(y)-, y) \, ({\bar s}(y)-)^{\gamma_1({\bar s}(y)-, y)}
+ C_2({\bar s}(y)-, y) \, ({\bar s}(y)-)^{\gamma_2({\bar s}(y)-, y)} \\
\notag
&= V({\bar s}(y), {\bar s}(y), y; \, b_*({\bar s}(y), y)) \equiv {\bar s}(y) - K 
\end{align}
for each ${\tilde y}_{2l}(s) < y \leq {\tilde y}_{2l-1}(s)$, $l \in \NN$, and $s > K$.
We can therefore conclude that the candidate value function admits the representation
\begin{align}
\label{4Vint31}
&V(x, s, y; {\bar s}(y), {\tilde y}_{2l-1}(s), {\tilde y}_{2l}(s)) \\
\notag
&= C_1(s, y; {\bar s}(y), {\tilde y}_{2l-1}(s), {\tilde y}_{2l}(s)) \, x^{\gamma_1(s, y)}
+ C_2(s, y; {\bar s}(y), {\tilde y}_{2l-1}(s), {\tilde y}_{2l}(s)) \, x^{\gamma_2(s, y)}
\end{align}
in the regions ${\tilde R}^3_{2l}$ given by (\ref{4tilR3_2l}), where
$C_i(s, y; {\bar s}(y), {\tilde y}_{2l-1}(s), {\tilde y}_{2l}(s))$, $i = 1, 2$, provide a
unique solution of the two-dimensional system of first-order linear partial differential equations in
(\ref{4B31c})-(\ref{4B31d}) with the boundary conditions of (\ref{4condcall1})-(\ref{4condcall2}),
for $l \in \NN$.
Finally, we observe that the candidate value function should be given by the condition
in the right-hand part of (\ref{4VD}) in the regions
\begin{equation}
\label{4hatR3_2k}
{\hat R}^3_{2k} = \{ (x,s,y) \in E^3 \, | \, {\hat y}_{2k}(s) < y \leq {\hat y}_{2k-1}(s) \}
\end{equation}
for $k \in \NN$, which belong to the stopping region $D_*$ in (\ref{4D}).

\vspace{3pt}


{\bf (ii) The case of put option.}
Let us finally consider the payoff function $G(x) = (L-x)^+$ in (\ref{4V5b}).
In this case, solving the system of equations in (\ref{4B31a}) and (\ref{4B31aa}), we conclude that the function in (\ref{4V0}) admits the representation
\begin{equation}
\label{4V32}
V(x, s, y; a_*(s, y)) = C_1(s, y; a_*(s, y)) \, x^{\gamma_1(s, y)} + C_2(s, y; a_*(s, y)) \, x^{\gamma_2(s, y)}
\end{equation}
for $0 < s-y \leq a_*(s, y) < x \le s$, with
\begin{equation}
\label{4Ci32}
C_i(s, y; a_*(s, y)) =
\frac{(\gamma_{3-i}(s,y) - 1) a_*(s, y) - \gamma_{3-i}(s, y) L}
{(\gamma_i(s, y) - \gamma_{3-i}(s,y)) a_*^{\gamma_i(s, y)}(s, y)}
\end{equation}
for all $0 < y < s$ and $i = 1, 2$.
Hence, assuming that the boundary function $a_*(s, y)$ is continuously
differentiable, we apply the condition of (\ref{4B31c}) for the functions
$C_i(s, y)= C_i(s, y; a_*(s, y))$, $i = 1, 2$, in (\ref{4Ci32})
to obtain that $a_*(s, y)$ solves the first-order nonlinear ordinary differential equation
\begin{align}
\label{4g'32}
\partial_s a(s, y) &= \sum_{i=1}^{2} \frac{((\gamma_{3-i}(s, y) - 1)
a(s, y) - \gamma_{3-i}(s, y)) a(s, y)}{(\gamma_i(s, y) - 1) (\gamma_{3-i}(s, y) - 1)
a(s, y) - \gamma_i(s, y) \gamma_{3-i}(s, y) \, L} \\
&\phantom{= \sum_{i=1}^{2} \;\:} \times \bigg( \frac{1}{\gamma_{3-i}(s,y) - \gamma_i(s,y)}
+ \frac{(s/a(s, y))^{\gamma_i(s, y)} \ln(s/a(s,y))}
{(s/a(s,y))^{\gamma_{i}(s,y)}-(s/a(s, y))^{\gamma_{3-i}(s,y)}}
\bigg) \, \partial_s \gamma_i(s,y) \nonumber
\end{align}
for $0 < y < s$, where the partial derivatives $\partial_s \gamma_{i}(s, y)$,
$i = 1, 2$, are given by (\ref{4gammas}) with (\ref{4phi}).

Since the functions $\delta(s, y)$ and $\sigma(s, y)$ are assumed to be
continuously differentiable and bounded, the limits $\delta(y+, y)$
and $\sigma(y+, y)$ exist for each $y > 0$. Then, the limits $\gamma_i(y+, y)$
can be identified with the functions $\beta_i(y)$, $i = 1, 2$, from Subsection 3.2 above,
and the function in (\ref{4V32}) should satisfy the property
$V(x, s, y; a_*(s, y)) \rightarrow V(x, y+, y; a_*(y+, y))$ as $s \downarrow y$,
for each $s-y \leq a_*(s, y) < x \le s$. Thus, we conclude that the equalities
\begin{equation}
\label{4Vyy}
V(x, y+, y; a_*(y+, y)) = U(x, y; a_*(y+, y)) \quad \text{and} \quad a_*(y+, y) = g_*(y)
\end{equation}
hold for $0 < a_*(y+,y) < x \leq y$ and $U(x, s; g_*(s))$ given by
(\ref{4U32a}) with $g_*(s)$ obtained in part (ii) of Subsection 3.2.
To see this, we observe that the candidate value function evaluated at $s \downarrow y$
in (\ref{4Vyy}) satisfies the normal reflection condition only at the diagonal
$d^3_3 = \{(x, s, y) \in \RR^3 \, | \, 0 < x = s = y \}$ of the plane $d^3_1$,
and thus, the function $a_*(y+, y) = g_*(y)$ is the maximal solution of the
equation in (\ref{4g'32a})
with the boundary condition $a_*(\infty, \infty) = g_*(\infty)$ of (\ref{4U32b})
as $y = s \rightarrow \infty$ and such that this solution stays strictly below the curve $x = L \wedge (r L/\delta(y))$.

For any $y > 0$ fixed, let us now consider the unique solution $a_*(s, y)$ of
(\ref{4g'32}) started at the value $a_*(y+, y) = g_*(y)$, given that
this solution stays strictly below the surface $x = L \wedge (r L/\delta(s, y))$. 
Then, we put ${\tilde s}_0(y) = y$ and define an increasing sequence $({\tilde s}_n(y))_{n \in \NN}$
such that the boundary $a_*(s, y)$ exits the region $E^3$ from the side of
the plane $d^3_1$ at the points $({\tilde s}_{2l-1}(y), {\tilde s}_{2l-1}(y), y)$ and
enters $E^3$ upwards at the points $({\tilde s}_{2l}(y), {\tilde s}_{2l}(y), y)$.
Namely, we define ${\tilde s}_{2l-1}(y) = \inf \{ s > {\tilde s}_{2l-2}(y) \, | \, a_*(s, y) > s \}$
and ${\tilde s}_{2l}(y) = \inf \{ s > {\tilde s}_{2l-1}(y) \, | \, a_*(s, y) \le s \}$,
$l \in \NN$, whenever they exist, and put ${\tilde s}_{2l-1}(y) = {\tilde s}_{2l}(y) = \infty$ otherwise, for $l \in \NN$.
Note that $y < {\tilde s}_{2l-1}(y) < {\tilde s}_{2l}(y) \leq L$, $l \in \NN$, by construction.
Moreover, we put ${\hat s}_0(y) = y$ and define an increasing sequence $({\hat s}_n(y))_{n \in \NN}$
such that ${\hat s}_{2k-1}(y) = \inf \{ s > {\hat s}_{2k-2}(y) \, | \, a_*(s, y) < s - y \}$
and ${\hat s}_{2k}(y) = \inf \{ s > {\hat s}_{2k-1}(y) \, | \, a_*(s, y) \ge s - y \}$,
$k \in \NN$, whenever they exist, and put ${\hat s}_{2k-1}(y) = {\hat s}_{2k}(y) = \infty$
otherwise. Note that $y \leq {\hat s}_{2k-2}(y) < {\hat s}_{2k-1}(y) < L+y$,
by construction, for $k = 1, \ldots, {\hat k}$, where
${\hat k} = \sup\{k \in \NN \, | \, {\hat s}_{2k-1}(y) - y < L\}$.
Therefore, the candidate value function admits the expression in (\ref{4V32})
in either the region
\begin{equation}
\label{4hatQ3_2k-2}
{\hat Q}^3_{2k-2} = \{ (x, s, y) \in E^3 \, | \, {\hat s}_{2k-2}(y) \leq s < \min_{l \in \NN}
\{ {\tilde s}_{2l-1}(y) | \, {\tilde s}_{2l-1}(y) > {\hat s}_{2k-2}(y) \} \wedge {\hat s}_{2k-1}(y) \}
\end{equation}
or
\begin{equation}
\label{4tilQ3_2l-2}
{\tilde Q}^3_{2l-2} = \{ (x, s, y) \in E^3 \, | \, {\tilde s}_{2l-2}(y) \leq y < \min_{k \in \NN}
\{ {\hat s}_{2k-1}(y) \, | \, {\hat s}_{2k-1}(y) > {\tilde s}_{2l-2}(y) \} \wedge {\tilde s}_{2l-1}(y) \}
\end{equation}
for $k = 1, \ldots, {\hat k}$, and $l \in \NN$,
and the boundary function $a_*(s, y)$ provides the unique solution of (\ref{4g'32}) started at the value
$a_*(y+,y) = g_*(y)$, given that this solution stays strictly below the surface $x = L \wedge (r L/\delta(s, y))$ (see Figure 2 below).

On the other hand, the candidate value function takes the form of (\ref{4V0}) with $C_{i}(s, y)$,
$i = 1, 2$, solving the linear system of first-order partial differential equations in (\ref{4B31c}) and (\ref{4B31d}),
in the regions
\begin{equation}
\label{4hatQ3_2k-1}
{\hat Q}^3_{2k-1} = \{ (x, s, y) \in E^3 \, | \, {\hat s}_{2k-1}(y) \leq s < {\hat s}_{2k}(y) \}
\end{equation}
for $k = 1, \ldots, {\hat k}$, which belong to $C'$ in (\ref{4C'}).
Note that, the process $(X, S, Y)$ can enter
${\hat Q}^3_{2k-1}$ in (\ref{4hatQ3_2k-1})
from one of the regions
${\hat Q}^3_{2k-2}$ in (\ref{4hatQ3_2k-2}) or ${\tilde Q}^3_{2l-2}$ in (\ref{4tilQ3_2l-2}),
for some $l \in \NN$,
only through the point $({\hat s}_{2k-1}(y), {\hat s}_{2k-1}(y), y)$ and can exit
${\hat Q}^3_{2k-1}$ passing to ${\hat Q}^3_{2k}$
only through the point
$({\hat s}_{2k}(y), {\hat s}_{2k}(y), y)$, by hitting the plane $d^3_1$
and increasing its second component $S$.
Thus,

\begin{picture}(160,110)
\put(20,15){\begin{picture}(120,95)

\put(0,0){\line(1,0){120}} \put(120,0){\line(0,1){95}}
\put(0,0){\line(0,1){95}} \put(0,95){\line(1,0){120}}   

\put(10,10){\vector(1,0){100}}
\put(10,10){\vector(0,1){80}}    



\put(9,35){\line(1,0){33}} 

\put(10.5,10){\line(5,4){95}} 

\put(6,34){${y}$}

\put(10,35){\line(5,4){65}} 

\put(48,40){\line(-1,0){19}}
\put(48,39.9){\line(-1,0){19}}
\put(48,40.1){\line(-1,0){19}}

\put(48,40){\line(5,4){12}}
\put(48,39.9){\line(5,4){12}}
\put(48,40.1){\line(5,4){12}}
\put(48,39.8){\line(5,4){12}}
\put(48,40.2){\line(5,4){12}}

\put(52,43){\line(-1,0){27}}
\put(52,42.9){\line(-1,0){27}}
\put(52,43.1){\line(-1,0){27}}


\put(56,46.5){\line(-1,0){15.5}}
\put(56,46.4){\line(-1,0){15.5}}
\put(56,46.6){\line(-1,0){15.5}}


\put(60,49.7){\line(-1,0){32}}
\put(60,49.5){\line(-1,0){32}}
\put(60,49.6){\line(-1,0){32}}


\qbezier[70](10,42)(37,64)(64,86)
\qbezier[40](10,41.5)(30,41.5)(50,41.5)
\qbezier[8](20,49.6)(24,49.6)(28,49.6)

\put(9,41.5){\line(1,0){2}} \put(6,41){$y'$}

\put(83,47){\vector(-1,1){25}} \put(73,44){$x = s-y$}

\put(34,40){${\bf .}$}
\put(34,39.9){${\bf .}$}	
\put(34.1,39.9){${\bf .}$}	
\put(33.9,39.9){${\bf .}$}
\put(34,39.8){${\bf .}$}	
\put(34.1,39.8){${\bf .}$} 	
\put(33.9,39.8){${\bf .}$}
\put(34,39.7){${\bf .}$}	
\put(34.1,39.7){${\bf .}$}	
\put(33.9,39.7){${\bf .}$}
\put(34,39.6){${\bf .}$}

\put(60,27){\vector(-2,1){24}} \put(61,26){$(x,s,y)$}

\put(6,88){$s$} \put(108,6){$x$}

\put(0,0){\qbezier(98,80)(80,84)(86,87)} 
\put(0,0){\qbezier(79,65)(111,75)(98,80)}
\put(0,0){\qbezier(35,55)(60,60)(79,65)}
\put(0,0){\qbezier(25,35)(18,50)(35,55)}

\qbezier[88](10,80)(54,80)(98,80)      
\qbezier[69](10,65)(44.5,65)(79,65)
\qbezier[25](10,55)(22.5,55)(35,55)
\qbezier[12](10,45.5)(16,45.5)(22,45.5)

\put(.5,79){${\tilde s}_{1}(y)$}
\put(.5,64){${\tilde s}_{2}(y)$}
\put(.5,54){${\hat s}_{1}(y)$}
\put(.5,45){${\hat s}_{2}(y)$}

\qbezier[25](25,35)(25,22.5)(25,10)
\put(25,9){\line(0,1){2}} \put(21,6){$g_*(y)$}

\put(105,9){\line(0,1){2}} \put(103,5){$L$}

\put(96,61){\vector(-1,2){4}} \put(96,58){$a_*(s,y)$}




\end{picture}}
\put(17,10){\small{{\bf Figure 2.} A computer drawing of the
state space of the process $(X,S,Y)$,}} \\
\put(17,5) {\small{for some $y$ fixed, which increases to
$y'$, and the boundary function $a_*(s,y)$.}}
\end{picture}

\noindent
the candidate value function should be continuous
at the points $({\hat s}_{2k-1}(y), {\hat s}_{2k-1}(y), y)$
and $({\hat s}_{2k}(y), {\hat s}_{2k}(y), y)$, that is expressed by the
equalities
\begin{align}
\label{4condput1}
&C_1({\hat s}_{2k-1}(y), y) \, ({\hat s}_{2k-1}(y))^{\gamma_1({\hat s}_{2k-1}(y), y)}
+ C_2({\hat s}_{2k-1}(y), y) \, ({\hat s}_{2k-1}(y))^{\gamma_2({\hat s}_{2k-1}(y), y)} \\
\notag
&= V({\hat s}_{2k-1}(y)-, {\hat s}_{2k-1}(y)-, y; a_*({\hat s}_{2k-1}(y)-, y)) \\
\label{4condput3}
&C_1({\hat s}_{2k}(y)-, y) \, ({\hat s}_{2k}(y)-)^{\gamma_1({\hat s}_{2k}(y)-, y)}
+ C_2({\hat s}_{2k}(y)-, y) \, ({\hat s}_{2k}(y)-)^{\gamma_2({\hat s}_{2k}(y)-, y)} \\
\notag
&= V({\hat s}_{2k}(y), {\hat s}_{2k}(y), y; a_*({\hat s}_{2k}(y), y))
\end{align}
for $y > 0$ and $k = 1, \ldots, {\hat k}-1$, where the right-hand sides are
given by (\ref{4V32})-(\ref{4Ci32}) with $a_*({\hat s}_{2k-1}(y)-, y) = ({\hat s}_{2k-1}(y) - y)-$
and $a_*({\hat s}_{2k}(y), y) = {\hat s}_{2k}(y) - y$, respectively.
Moreover, in the region
${\hat Q}^3_{2{\hat k}-1}$ we have ${\hat s}_{2{\hat k}}(y) = \infty$
and the condition of (\ref{4condput3}), for $k={\hat k}$, changes its form to
$C_1(\infty, y) = 0$ for $y > 0$, since otherwise $V(x, \infty, y) \rightarrow \pm
\infty$ as $x \uparrow \infty$,
that must be excluded due to the fact that the value function in (\ref{4V5b}) is bounded at infinity,
while the condition of (\ref{4condput1}) holds for $k={\hat k}$ as well.

In addition, the process $(X, S, Y)$ can exit
${\hat Q}^3_{2k-1}$ in (\ref{4hatQ3_2k-1}) passing to the stopping region $D_*$ in (\ref{4D}),
only through the point $(s-{\bar y}(s), s, {\bar y}(s))$, by hitting the plane $d^3_2$,
so that increasing its third
component $Y$ until it reaches the value ${\bar y}(s) = \inf \{ z > y \, | \, a_*(s, z) \geq s - z \}$.
Observe that the boundary $a_*(s, z)$ provides the unique solution of the equation in (\ref{4g'32}) with
the starting value $a_*(z+, z) = g_*(z)$ from (\ref{4g'32a}), for each $z < {\bar y}(s)$,
given that this solution stays strictly below the surface $x = L \wedge (r L/\delta(s, z))$. 
Then, the candidate value function should be 
continuous at the point $(s-{\bar y}(s), s, {\bar y}(s))$,
that is expressed by the equality
\begin{align}
\label{4condput2}
&C_1(s, {\bar y}(s)-) \, ((s-{\bar y}(s))+)^{\gamma_1(s, {\bar y}(s)-)} +
C_2(s, {\bar y}(s)-) \, ((s-{\bar y}(s))+)^{\gamma_2(s, {\bar y}(s)-)}
\\ \notag
&= V(s-{\bar y}(s), s, {\bar y}(s); a_*(s, {\bar y}(s))) \equiv L - (s - {\bar y}(s)) 
\end{align}
for each ${\hat s}_{2k-1}(y) \leq s < {\hat s}_{2k}(y)$, $k = 1, \ldots, {\hat k}$, and $y > 0$.
We can therefore conclude that the candidate value function admits the representation
\begin{align}
\label{4Vint32}
&V(x, s, y; {\hat s}_{2k-1}(y), {\hat s}_{2k}(y), {\bar y}(s)) \\
\notag
&= C_1(s, y; {\hat s}_{2k-1}(y), {\hat s}_{2k}(y), {\bar y}(s)) \, x^{\gamma_1(s, y)} +
C_2(s, y; {\hat s}_{2k-1}(y), {\hat s}_{2k}(y), {\bar y}(s)) \, x^{\gamma_2(s, y)}
\end{align}
in the regions ${\hat Q}^3_{2k-1}$ in (\ref{4hatQ3_2k-1}),
where $C_i(s, y; {\hat s}_{2k-1}(y), {\hat s}_{2k}(y), {\bar y}(s))$, $i = 1, 2$,
provide a unique solution of the two-dimensional system of linear partial
differential equations in (\ref{4B31c})-(\ref{4B31d}) with the boundary conditions
(\ref{4condput1})-(\ref{4condput2}), for $k = 1, \ldots, {\hat k}$.
Finally, we note that the candidate value function should be given by the condition
in the left-hand part of (\ref{4VD}) in the regions
\begin{equation}
\label{4tilQ3_2l-1}
{\tilde Q}^3_{2l-1} = \{ (x, s, y) \in E^3 \, | \, {\tilde s}_{2l-1}(y) \leq s < {\tilde s}_{2l}(y) \}
\end{equation}
for $l \in \NN$, which belong to the stopping region $D_*$ from (\ref{4D}).


\section{\dot Main results and proof}
  \setcounter{equation}{0}

\label{ch4.3}

     In this section, taking into account the facts proved above,
     we formulate and prove the main results of the paper. We recall
     that the process $(X, S, Y)$ is defined in (\ref{4X4})-(\ref{4S4}).


\begin{theorem} \label{4call}
In the perpetual American call option case with the payoff function $G(x) = (x-K)^+$,
the value function of the optimal stopping problem (\ref{4V5b}) has the expression
\begin{equation}
\label{4V*4c}
V_*(x, s, y)=
\begin{cases}
V(x, s, y; b_*(s, y)), \;\; 
&\text{if} \;\; s-y \leq x < b_*(s, y) \leq s \\
V(x, s, y; {\bar s}(y), {\tilde y}_{2l-1}(s), {\tilde y}_{2l}(s)), \;\; 
&\text{if} \;\; s-y \leq x \leq s < b_*(s, y) \\
x - K, \;\; 
&\text{if} \;\; b_*(s, y) \leq x \leq s
\end{cases}
\end{equation}
and the optimal stopping time is given by the right-hand part of (\ref{4tau}),
where the functions $V(x, s, y; b_*(s, y))$ and $V(x, s, y; {\bar s}(y),
{\tilde y}_{2l-1}(s), {\tilde y}_{2l}(s))$ as well as the boundary $b_*(s, y)$
are specified as follows:

(i) in the particular case $\delta(s, y) = \delta(s)$ and $\sigma(s, y) = \sigma(s)$, the function
$V(x, s, y; b_*(s, y)) = U(x, s; h_*(s))$ and the boundary $b_*(s, y) = h_*(s)$ are given
by (\ref{4U33a}), for $(x, s) \in {\tilde R}^2_{2l-1}$ defined in (\ref{4tilR2_2l-1}),
and $V(x, s, y; {\bar s}(y),
{\tilde y}_{2l-1}(s), {\tilde y}_{2l}(s)) = U(x, s; {\tilde s}_{2l-1})$ is given by
(\ref{4Uint2}), whenever $(x, s) \in {\tilde R}^2_{2l}$ defined in (\ref{4tilR2_2l}),
for $l \in \NN$;

(ii) in the general case for $\delta(s, y)$ and $\sigma(s, y)$, the function $V(x, s, y; b_*(s, y))$
is given by (\ref{4V33})-(\ref{4Ci33}) and the boundary $b_*(s, y)$ provides the unique solution of
the equation in (\ref{4g'33}) started at the value $b_*(s, s-) = h_*(s)$ from (\ref{4U33a}),
given that $b_*(s, y) > K \vee (r K/\delta(s, y))$ holds
for $(x, s, y) \in {\tilde R}^3_{2l-1} \cup {\hat R}^3_{2k-1}$ defined in
(\ref{4tilR3_2l-1}) and (\ref{4hatR3_2k-1}), respectively,
and $V(x, s, y; {\bar s}(y), {\tilde y}_{2l-1}(s),$  ${\tilde y}_{2l}(s))$
is given by (\ref{4Vint31}), whenever $(x, s, y) \in {\tilde R}^3_{2l}$
defined in (\ref{4tilR3_2l}), with $C_i(s, y; {\bar s}(y), {\tilde y}_{2l-1}(s),$
${\tilde y}_{2l}(s))$, $i = 1, 2$, solving the system of equations in
(\ref{4B31c})-(\ref{4B31d}) and satisfying the conditions of
(\ref{4condcall1})-(\ref{4condcall2}), for $k, l \in \NN$,
where (\ref{4condcall3}) changes its form to $C_2(s, s-) = 0$,
for the case $l = 1$, if $b_*(s,s-) = h_*(s) > s$ holds.
\end{theorem}


\begin{theorem} \label{4put}
In the perpetual American put option case with the payoff function $G(x) = (L - x)^+$,
the value function of the optimal stopping problem (\ref{4V5b}) has the expression
\begin{equation}
\label{4V*4p}
V_*(x, s, y)=
\begin{cases}
V(x, s, y; a_*(s, y)), \;\; 
&\text{if} \;\; s-y \leq a_*(s, y) < x \leq s\\
V(x, s, y; {\hat s}_{2k-1}(y), {\hat s}_{2k}(y), {\bar y}(s)), \;\; 
&\text{if} \;\; a_*(s, y) < s-y \leq x \leq s\\
L - x, \;\; 
&\text{if} \;\; s-y \leq x \le a_*(s, y)
\end{cases}
\end{equation}
and the optimal stopping time is given by the left-hand part of (\ref{4tau}), where
the functions $V(x, s, y; a_*(s, y))$ and $V(x, s, y; {\tilde s}_{2l-1}(y), {\tilde s}_{2l}(y),
{\bar y}(s))$ as well as the boundary $a_*(s, y)$ are specified as follows:

(i) in the particular case $\delta(s, y) = \delta(s)$ and $\sigma(s, y) = \sigma(s)$, the function
$V(x, s, y; a_*(s, y)) = U(x, s; g_*(s))$ is given by (\ref{4U32a})-(\ref{4Di32a}) and the
boundary $a_*(s, y) = g_*(s)$ provides the maximal solution of the equation in (\ref{4g'32a}) started
at $g_*(\infty)$ from (\ref{4U32b}), such that $g_*(s) < L \wedge (r L/\delta(s))$ holds
for $(x, s) \in {\hat Q}^2_{2k-1}$ defined in (\ref{4hatR2_2k-1}) and $k \in \NN$;

(ii) in the general case for $\delta(s, y)$ and $\sigma(s, y)$, the function $V(x, s, y; a_*(s, y))$
is given by (\ref{4V32})-(\ref{4Ci32}) and the boundary $a_*(s, y)$ provides the unique solution
of the equation in (\ref{4g'32}) started from the value
$a_*(y+, y) = g_*(y)$ from part (i) above, given that $a_*(s, y) < L \wedge (r L/\delta(s, y)$ holds
for $(x, s, y) \in {\hat Q}^3_{2k-2} \cup {\tilde Q}^3_{2l-2}$ defined in
(\ref{4hatQ3_2k-2}) and (\ref{4tilQ3_2l-2}), respectively,
and $V(x, s, y; {\hat s}_{2k-1}(y), {\hat s}_{2k}(y), {\bar y}(s))$
is given by (\ref{4Vint32}), whenever $(x, s, y) \in {\hat Q}^3_{2k-1}$ defined in
(\ref{4hatQ3_2k-1}), with
$C_i(s, y; {\hat s}_{2k-1}(y), {\hat s}_{2k}(y), {\bar y}(s))$, $i = 1, 2$,
solving the system of equations in (\ref{4B31c})-(\ref{4B31d}) and satisfying
the conditions of (\ref{4condput1})-(\ref{4condput2}),
$k = 1, \ldots, {\hat k}$, and $l \in \NN$,
where (\ref{4condput3}) changes its form to $C_1(\infty, y) = 0$, for the case $k = {\hat k}$.
\end{theorem}


Since all the parts of the assertions formulated above are proved using similar arguments, we only
give a proof for the three-dimensional optimal stopping problem related to the perpetual American
put option in part (ii) of Theorem \ref{4put}, which represents the most complicated and informative case.


\begin{proof} {\bf of part (ii) of Theorem \ref{4put}.}
In order to verify the assertion stated above, it remains to show
        that the function defined in (\ref{4V*4p}) coincides with the value function in
        (\ref{4V5b}) and that the stopping time $\tau_*$ in the left-hand part of
(\ref{4tau}) is optimal with
        the boundary $a_*(s, y)$ specified above. For this,
let $a(s, y)$ be the unique solution of (\ref{4g'32}) starting from the value $a(y+,y) = g(y)$, being any solution of (\ref{4g'32a}) starting from $a_*(\infty, \infty) = g_*(\infty)$ in (\ref{4U32b}) and satisfying $g(y) < L \wedge (r L/\delta(y))$ for all $y > 0$. Let us also denote by $V_{a}(x, s, y)$ the right-hand side of the expression
in (\ref{4V*4p}) associated with this $a(s, y)$.
		It then follows using straightforward calculations and the assumptions
		presented above that the function $V_{a}(x, s, y)$ solves the left-hand system of
        (\ref{4LV})-(\ref{4VD}),
        while the normal-reflection and smooth-fit conditions are satisfied
        in (\ref{4NRS}) and the left-hand part of (\ref{4NRY}).
        Hence, taking into account the fact that the function $V_{a}(x, s, y)$
        is $C^{2,1,1}$ and the boundary $a(s, y)$ is assumed to be
        continuously differentiable for all $0 < y < s$, by applying the change-of-variable
        formula from \cite[Theorem~3.1]{Pe1a} to $e^{- r t} \, V_{a}(X_t, S_t, Y_t)$, we obtain
        \begin{align}
        \label{4rho4c}
        &e^{- r t} \, V_{a}(X_t, S_t, Y_t) = V_{a}(x, s, y) + M_t \\
        \notag
        &+ \int_0^t e^{- r u} \,
        (\LL V_{a} - r V_{a}) (X_u, S_u, Y_u) \, I(X_u \not= S_{u} - Y_u, X_u \not= a(S_u, Y_u), X_u \not= S_u) \, du \\
        \notag
        &+ \int_0^t e^{- r u} \, \partial_s V_{a}(X_u, S_u, Y_u) \, I(X_u = S_u) \, dS_u 
        + \int_0^t e^{- r u} \, \partial_y V_{a}(X_u, S_u, Y_u) \, I(X_u = S_u - Y_u) \, dY_u
        \end{align}
        where $I(\cdot)$ denotes the indicator function and the process $M = (M_t)_{t \ge 0}$ given by
        \begin{equation}
        \label{4N5}
        M_t = \int_0^t e^{- r u} \, \partial_x V_{a}(X_u, S_u, Y_u) \,
        I(X_u \neq S_{u} - Y_u, X_u \neq S_{u})
        \, \sigma(S_u, Y_u) \, X_{u} \, dB_u
        \end{equation}
        is a square integrable martingale under $P_{x, s, y}$.
        Note that, since the time spent by the process $X$ at the boundary surface
        $\{ (x, s, y) \in E^3 \, | \, x = a(s, y) \}$
as well as at the planes $d^3_1$ and $d^3_2$,
        is of Lebesgue measure zero,
        the indicators in the second line of the formula (\ref{4rho4c}) as well as
        in the formula (\ref{4N5}) can be ignored.
        Moreover, since the process $S$ increases only at the plane $d^3_1$
        and the process $Y$ increases only at the plane $d^3_2$,
        the indicators in the third and fourth line of (\ref{4rho4c}) can also be set equal to one.

        By using straightforward calculations and the arguments from the previous section,
        it is verified that $(\LL V_{a} - r V_{a})(x, s, y) \le 0$ for all $(x, s, y) \in E^3$
        such that $x \neq a(s, y)$,
$x \neq s - y$, and $x \neq s$.
        Moreover, it is shown by means of standard arguments that the property
on the left-hand part of (\ref{4VC})
        also holds, which together with the left-hand parts of (\ref{4CF})-(\ref{4VD}) imply that the inequality
        $V_{a}(x, s, y) \ge (L - x)^+$
is satisfied for all $(x, s, y) \in E^3$.
        It therefore follows from the expression (\ref{4rho4c}) that the inequalities
        \begin{equation}
        \label{4rho4e}
        e^{- r \tau} \, 
(L - X_{\tau})^+
\le e^{- r \tau} \, V_{a}(X_{\tau}, S_{\tau}, Y_{\tau}) \le V_{a}(x, s, y) + M_{\tau}
        \end{equation}
        hold for any finite stopping time $\tau$ with respect to the natural filtration of $X$.

        Taking the expectation with respect to $P_{x, s, y}$ in (\ref{4rho4e}), by means of the optional sampling theorem
        (see, e.g. 
		\cite[Chapter~I, Theorem~3.22]{KS}), we get
        \begin{align}
        \label{4rho4e2}
        E_{x, s, y} \big[ e^{- r (\tau \wedge t)} \,
(L - X_{\tau \wedge t})^+
\big] &\le
        E_{x, s, y} \big[ e^{- r (\tau \wedge t)}
        \, V_{a}(X_{\tau \wedge t}, S_{\tau \wedge t}, Y_{\tau \wedge t}) \big] \\
        \notag
        &\le V_{a}(x, s, y) + E_{x, s, y} \, M_{\tau \wedge t} = V_{a}(x, s, y)
        \end{align}
        for all $(x, s, y) \in E^3$.
        Hence, letting $t$ go to infinity and using Fatou's lemma,
        we obtain that for any finite stopping time $\tau$ the inequalities
        \begin{equation}
        \label{4rho4e3}
        E_{x, s, y} \big[ e^{- r \tau} \, 
(L - X_{\tau})^+
\big] \le E_{x, s, y} \big[ e^{- r \tau} \, V_{a}(X_{\tau}, S_{\tau}, Y_{\tau}) \big] \le V_{a}(x, s, y)
        \end{equation}
        are satisfied for all $(x, s, y) \in E^3$.
Taking first the supremum over all stopping times $\tau$ and then the infimum over all $a$, we conclude that
        \begin{equation}
        \label{4rho4e32}
E_{x, s, y} \big[ e^{- r \tau_*} \, (L - X_{\tau_*})^+ \big]
\le \inf_{a} V_{a}(x, s, y) = V_{a_*}(x, s, y)
        \end{equation}
where $a_*(s,y)$ is the unique solution of (\ref{4g'32}) starting from the value $a_*(y+,y) = g_*(y)$,
which is the maximal solution to (\ref{4g'32a}) being started at $a_*(\infty, \infty) = g_*(\infty)$
in (\ref{4U32b}) and staying strictly below the curve $x = L \wedge (r L/\delta(y))$.
Using the fact that $V_{a}(x, s, y)$ is decreasing in the function $a < L \wedge (r L/\delta)$, we see that the infimum in (\ref{4rho4e32}) is attained over any sequence of solutions $(a_n(s, y))_{n \in \NN}$ to (\ref{4g'32}) started at the values $a_n(y+, y) = g_n(y)$, solving (\ref{4g'32a}) and such that $g_n(y) \uparrow g_*(y)$, and thus, $a_n(s,y) \uparrow a_*(s,y)$
as $n \rightarrow \infty$. Since the inequalities in (\ref{4rho4e3}) hold also for $a_*(s, y)$, we see that (\ref{4rho4e32})
holds for $a_*(s,y)$ and $(x,s,y) \in E^3$ as well. Note that $V_{a}(x, s, y)$ in (\ref{4rho4e2}) is superharmonic
for the Markov process $(X,S,Y)$ on $E^3$.
Taking into account the fact that $V_{a}(x, s, y)$ is decreasing in $a < L \wedge (r L/\delta)$ and
that the inequality $V_{a}(x, s, y) \geq (L-x)^+$ holds for all $(x,s,y) \in E^3$,
we observe that the selection of the maximal solution $a_*(s,y)$,
which stays strictly below the surface $x = L \wedge (r L/\delta(s, y))$, whenever such a choice exists,
is equivalent to invoking the superharmonic characterization of
the value function (smaller superharmonic function dominating
the payoff function, see \cite{Pmax} or \cite[Chapter~1]{PSbook}).

In order to prove the fact that $a_*(s,y)$ is optimal on $E^3$, we consider the sequence of stopping times $\tau_n$ defined as in the left-hand part of (\ref{4tau}) with $a_n(s,y)$ instead of $a_*(s,y)$, where $a_n(s,y)$ is the unique solution of (\ref{4g'32}) started from the value $a_n(y+,y) = g_n(y)$ which solves (\ref{4g'32a}) and starts at $a_*(\infty, \infty) = g_*(\infty)$ in (\ref{4U32b}), and such that $g_n(y_n) = L \wedge (r L/\delta(y_n))$, for some $y_n \downarrow 0$ as $n \rightarrow \infty$.
By virtue of the fact that the function $V_{a_n}(x, s, y)$ from the right-hand side of the expression in (\ref{4V*4p}) associated with this $a_n(s, y)$, satisfies the left-hand system of
        (\ref{4LV})-(\ref{4VC}) with (\ref{4NRS}) and taking into account the structure
        of $\tau_n$ given by the left-hand part of (\ref{4tau}) with $a_n(s,y)$ instead of $a_*(s,y)$, it follows from the equivalent expression of (\ref{4rho4c})
        that the equalities
        \begin{align}
        \label{4rho4g}
        e^{- r (\tau_n \wedge t)} \, 
(L - X_{\tau_n \wedge t})^+
        = e^{- r (\tau_n \wedge t)} \, V_{a_n}(X_{\tau_n \wedge t}, S_{\tau_n \wedge t}, Y_{\tau_n \wedge t})
        = V_{a_n}(x, s, y) + M_{\tau_n \wedge t}
        \end{align}
        hold for all $(x, s, y) \in E^3$.
        Observe that
        $\tau_n \uparrow \tau_*$ ($P_{x, s, y}$-a.s.)
and the variable $e^{- r \tau_*} (L - X_{\tau_*})^+$
        is bounded on the set $\{ \tau_* = \infty \}$. Taking into account the fact that the boundary $a_*(s, y)$
is bounded, it is easily seen that the property $P_{x, s, y}(\tau_* < \infty) = 1$ holds
        for all $(x, s, y) \in E^3$.
        Hence, letting $t$ and $n$ go to infinity and using the conditions on the left-hand part of
        (\ref{4CF}) and in (\ref{4NRS}), as well as the fact that $\tau_n \uparrow \tau_*$ ($P_{x, s, y}$-a.s.),
        we can apply the Lebesgue dominated convergence theorem for (\ref{4rho4g}) to obtain the equality
        \begin{equation}
        \label{4rho4i}
        E_{x, s, y} \big[ e^{- r \tau_*} \, 
(L - X_{\tau_*})^+
\big] = V_{a_*}(x, s, y)
        \end{equation}
        for all $(x, s, y) \in E^3$, which together with (\ref{4rho4e32})
        directly implies the desired assertion. 
\end{proof}

{\bf Acknowledgments.}
The second author gratefully acknowledges the scholarship of the Alexander Onassis Public Benefit Foundation for
his doctoral studies at the London School of Economics and Political Science.

\end{document}